\documentclass[a4paper]{article}
\title{A 64-dimensional two-distance counterexample to Borsuk's conjecture}
\author{Thomas Jenrich}
\date{2014-08-20}
\begin{document}
\maketitle

\section{Abstract and introduction}

In \cite{Bor} (1933) Karol Borsuk asked whether each bounded set in the
$n$-dimensional Euclidean space can be divided into $n$+1 parts of smaller diameter.
The diameter of a set is defined as the supremum (least upper bound) of the
distances of contained points. Implicitly, the whole set is assumed to
contain at least two points.

The hypothesis that the answer to that question is positive became famous
under the name \emph{Borsuk's conjecture}. Beginning with Jeff Kahn and Gil
Kalai, from 1993 to 2003 several authors have proved that in certain (almost
all) high dimensions such a division is not generally possible.

In \cite{Bon} (2013) and finally \cite{Bon2} Andriy V. Bondarenko constructed
a 65-dimensional two-distance set of 416 vectors that cannot be divided into
less than 84 parts of smaller diameter. That was not just the first known
two-distance counterexample to Borsuk's conjecture but also a considerable
reduction of the lowest known dimension the conjecture fails in in general.

This article presents a 64-dimensional subset of the vector set mentioned
above that cannot be divided into less than 71 (by \cite{Bon2} 72) parts of
smaller diameter, that way delivering a two-distance counterexample to
Borsuk's conjecture in dimension 64.

The contained proof relies on the results of some (combinatorial)
calculations. The additionally (in the source package) provided small computer
program G24CHK needs about one second for that task on a 1 GHz Intel PIII.

Meanwhile a short paper by this author and Andries E. Brouwer that follows
the principal idea of this article but avoids the extensive computational
part has been submitted to The Electronic Journal of Combinatorics.

\section{Euclidean representations of strongly regular graphs}

The following is mainly an essence of the content of the corresponding
section in \cite{Bon}, giving just the needed basic facts.
Detailed information can be found e.g. in \cite{BroHae}.

Saying $G$ is a srg$(v,k,\lambda,\mu)$, where srg abbreviates
``strongly regular graph'', means that there is a set $V$ of $v$ elements
(vertices) such that $G \subset V \times V$, for all $i \in V$
$\lnot G(i,i) \land |\{j: G(i,j)\}|=k$, and for all different $i,j \in V$
$$G(i,j)=G(j,i)
$$
and
$$
|\{ p \in V : G(p,i) \land G(p,j)\}| = \left\{\begin{array}{ll}
\lambda & \mbox{if $G(i,j)$}\\
\mu & \mbox{otherwise}\end{array}\right.
$$

The adjacency matrix $A$ of $G$ has exactly 3 different eigenvalues:
$k$ of multiplicity $1$, one positive eigenvalue of multiplicity
$$
f=\frac 12
\left(v-1-\frac{2k+(v-1)(\lambda-\mu)}{\sqrt{(\lambda-\mu)^2+4(k-\mu)}}\right),
$$
and one negative eigenvalue
$$
s=\frac 12\left(\lambda-\mu-\sqrt{(\lambda-\mu)^2+4(k-\mu)}\right)
$$

In the remaining part of this article we use these notations:
$I$ is the identity matrix of size $v$, $y$ is $A-sI$, $y_i$,
where $i \in V$, are the columns of $y$, and $y_{i,j}$, where
$i,j \in V$, are the entries of $y$.

The given properties of the eigenvalues imply $\dim \{y_i : i \in V\} \leq f$.

\section{The $G_2(4)$ graph}

Let $G$ be isomorphic to a graph known as $G_2(4)$ graph, a srg(416,100,36,20).
The value of $s$ is $-4$, thus for $i, j \in V$
$$
y_{i,j} = \left\{\begin{array}{ll}
4 & \mbox{if $i=j$}\\
1 & \mbox{if $G(i,j)$}\\
0 & \mbox{otherwise}\end{array}\right.
$$

For $i\in V$, $y_i$ consists of one 4 (at position $i$), 100 1's, and
315 0's.
For different $i,j\in V$
$$
\| y_i - y_j \| ^ 2 = \left\{\begin{array}{rl}
2\times{(100-36-1+(4-1)^2)} = 144  & \mbox{if $G(i,j)$}\\
2\times{(100-20+4^2)} = 192 & \mbox{otherwise}\end{array}\right.
$$
The value of $f$ is 65 and therefore $\dim \{y_i: i\in V\} \le 65$.

For each $W \subseteq V$: $\{y_i: i\in W\}$ has a smaller diameter than $\{y_i: i\in V\}$ iff all members of $W$ are pairwise adjacent.

By \cite{Bon} the $G_2(4)$ graph does not contain a clique (subgraph of
pairwise adjacent vertices) with more than 5 vertices, and because there
are 416 vertices (and corresponding vectors), $\{y_i: i\in V\}$ cannot be
divided into less than 84 subsets of smaller diameter.

\section{A construction of the $G_2(4)$ graph}

The following algorithm is an extract of the description
given in \cite{Bro1} and \cite{Bro2}, there originated to \cite{HoKiNa},
meanwhile arrived in the joint paper \cite{BrHoKiNa}.
See also \cite{Soi} and \cite{CrnMik}. Notice: The points mentioned here
are strictly to distinct from the points in the $n$-dimensional Euclidean
space.

``Consider the projective plane PG(2,16) provided with a nondegenerate
Hermitean form. It has 273 points, 65 isotropic and 208 nonisotropic.
There are 416 = 208 $\cdot$ 12 $\cdot$ 1/6  orthogonal bases. These are the
vertices of $G$. [...] Associated with a basis \{a,b,c\} is the triangle
consisting of the 15 isotropic points on the three lines $ab$, $ac$, and $bc$.
[...] $G$ can be described as the graph on the 416 triangles, adjacent when
they have 3 points in common.''

\section{A division of the vertices of the $G_2(4)$ graph}

The construction given in the previous section assigns a set of 15 different
isotropic points to each of the 416 vertices/triangles. Assume that integer
numbers from 1 to 65 are given to the isotropic points and consider the
sets of the numbers of the isotropic points, here shortly named iso-sets.
Thus two vertices are adjacent if the cardinality of the intersection of
their iso-sets is 3.
Now we start to divide the set of the 416 vertices:
First we define $B$ to be the set of all vertices in $V$ containing the
number 1 in their iso-set and $C$ to be the set of the remaining vertices.
The edges of $G$ connecting two vertices in $B$ constitute a subgraph of $G$.
We divide $B$ into non-empty subsets $B_h$, $h$ positive integer, such that
two different vertices are in the same subset iff $G$ contains a path between
them not leaving $B$.

\section{The computer program G24CHK}

The computer program G24CHK, provided in the source package of this article
and described in more detail in the last sections before the references,
implements the construction of the iso-sets of the members of $V$ and its
just defined subsets. Additionally, it checks that $B$ consists
of 3 pairwise disjunctive subsets ($B_1$, $B_2$, and $B_3$), each containing
exactly 32 vertices. Consequently, $|B| = 96 \land |C| = 416-96=320$.

Finally, G24CHK is able to check this claim:

(1) $\forall i \in V, h \in \{1,2,3\}$
$$
| \{j \in B_h : G(i,j) \} | = \left\{\begin{array}{ll}
20 & \mbox{if $i \in B_h$}\\
0 & \mbox{if $i \in B \setminus {B_h}$}\\
8 & \mbox{otherwise}\end{array}\right.
$$

The second case is already implied by the construction of the $B_h$.
The number and cardinalities of the unconnected subsets of $B$ and the first
two cases in (1) correspond with this slightly cut quotation from \cite{Bro2}, section
\emph{Subgraphs}:

``b) \emph{Three copies of the 2-coclique extension of the Clebsch graph.}

The 65520 nonedges of $G$ [...] fall into 1365 sets of 48, where each set
induces a subgraph of size 96 that is the disjoint union of three copies of
the 2-coclique extension of the Clebsch graph (the halved 5-cube), and
the 48 nonedges are the 2-cocliques.''

So actually just for the proof of the third case in (1) the actual
construction and counting are needed here.

However, in the following two sections the claims concerning $G$ and the
subsets of $V$ are taken as facts.

\section{The 64-dimensional counterexample}

Recall that for $i, j \in V$
$$
y_{i,j} = \left\{\begin{array}{ll}
4 & \mbox{if $i=j$}\\
1 & \mbox{if $G(i,j)$}\\
0 & \mbox{otherwise}\end{array}\right.
$$
Therefore statement (1) implies

(2)
$\forall i \in V, h \in \{1,2,3\}$

$$
\sum_{j \in B_h} y_{i,j} = \left\{\begin{array}{rl}
20 \times 1 + 1 \times 4 = 24 & \mbox{if $i \in B_h$}\\
0 & \mbox{if $i \in B \setminus {B_h}$}\\
8 \times 1 = 8 & \mbox{otherwise}\end{array}\right.
$$

We define a vector $p$ in the same 416-dimensional space as the
$y_i$, $i \in V$, consisting of entries $p_j$, $j \in V$:
$$
p_j = \left\{\begin{array}{rl}
1 & \mbox{if $j \in B_2$}\\
-1 & \mbox{if $j \in B_3$}\\
0 & \mbox{otherwise}.\end{array}\right.
$$

Combined with (2) this definition implies $\forall i \in V$
$$
\langle p, y_i \rangle = \sum_{j\in B_2} y_{i,j} - \sum_{j\in B_3} y_{i,j} =\\ \left\{\begin{array}{rl}
0 - 0  = 0 & \mbox{if $i \in B_1$}\\
24 - 0 = 24 & \mbox{if $i \in B_2$}\\
0 - 24 = -24 & \mbox{if $i \in B_3$}\\
8 - 8 = 0 & \mbox{otherwise}\end{array}\right.
$$

The essence is that $p$ is orthogonal to \{$y_i : i \in C \cup B_1\}$,
but not to \{$y_i : i \in V\}$.
Therefore $ \dim \{y_i : i \in C\cup B_1\} \le \dim \{y_i, i \in V\} - 1 $.
Because $\dim \{y_i : i \in V\} \le 65$ is known, we have $\dim \{y_i : i \in C\cup B_1\} \le 64$.

The proofs are not included here, but one can show (e.g. by actual vector
calculations) that these statements stay valid if we replace $\le$ by $=$.

Because \{$y_i : i \in C\cup B_1\}$ contains 352 vectors and a subset of
smaller diameter contains at most 5 vectors, a division into less than 71
parts of smaller diameter is impossible.
Because $71 > 64+1 $, the answer to Borsuk's question for $n=64$ is
negative.

In \cite{Bon2} A. Bondarenko states in a remark that a computer check had
shown that the number of parts of smaller diameter is at least 72.

\section{A 63-dimensional almost-counterexample}

We define a vector $q$ in the same 416-dimensional space as the
$y_i$, $i \in V$, consisting of entries $q_j$, $j \in V$:
$$
q_j = \left\{\begin{array}{rl}
2 & \mbox{if $j \in B_1$}\\
-1 & \mbox{if $j \in B_2 \cup B_3$}\\
0 & \mbox{otherwise}.\end{array}\right.
$$

Observe that $\langle p, q \rangle =  0 - 32 + 32 = 0$.

Combined with (2) the definition of $q$ implies $\forall i \in V$
$$
\langle q, y_i \rangle = 2 \times \sum_{j\in B_1} y_{i,j} - \\
 \sum_{j\in B_2} y_{i,j}-\sum_{j\in B_3} y_{i,j} =\\
\left\{\begin{array}{rl}
48 - 0 - 0  = 48 & \mbox{if $i \in B_1$}\\
0 - 24 - 0 = -24 & \mbox{if $i \in B_2$}\\
0 - 0 - 24 = -24 & \mbox{if $i \in B_3$}\\
16 - 8 - 8 = 0 & \mbox{otherwise}\end{array}\right.
$$

The essence is that $q$ is orthogonal to \{$y_i : i \in C \}$,
but not to \{$y_i : i \in C \cup B_1 \}$.
Therefore $ \dim \{y_i : i \in C\} \le \dim \{y_i : i \in C\cup B_1\} - 1 $.
Because $\dim \{y_i : i \in C\cup B_1\} \le 64$ has been derived, we have
$\dim \{y_i : i \in C\} \le 63$.

Again, these statements stay valid if we replace $\le$ by $=$.

But \{$y_i : i \in C\}$ contains just 320 vectors.
The proof is not included here, but one can divide $C$ into 64 5-cliques
and therefore \{$y_i : i \in C\}$ into 64 parts of smaller diameter.

As checked computationally, one such partition uniquely fulfills the
additional requirement that for each of the 64 5-cliques the (common)
intersection of the iso-sets of all five contained vertices has size 3.
Because of this property, that partition is relatively easy to find.

\section{Mathematical foundations of the implementation in G24CHK}

The $G_2(4)$ graph construction algorithm given above is very concise.
Sufficient information on projective planes can be found in \cite{Che}.
The implementation of the algorithm uses elements of the three-dimensional
vector space over the finite field GF(16) to constitute and represent
the objects in the projective plane PG(2,16). Enough information on finite
fields can be found in \cite{Cam}, using just GF(16) as example.
Unfortunately, the calculation of the list of the polynomials representing the
elements of GF(16) contains an error, resulting especially in the
impossible relation $x^{11} = x^{13}$.

The applied Hermitean form takes three-dimensional vectors $a$ and $b$ over GF(16)
and returns \mbox{$a_1 \overline {b_3} + a_2 \overline {b_2} + a_3 \overline {b_1}$},
where addition, multiplication, and conjugation operate in GF(16).

Finally, \cite{CibLip} contains enough information on isotropic objects.

\section{Compiling and executing G24CHK}

To compile the provided source code file G24CHK.PAS you will need a
PASCAL compiler compatible with Turbo Pascal 4.0. Even Turbo Pascal versions
3.02 and 1.0 are sufficient after removing the right half of the first
occurrence of the sequence (**), this way transforming the source up to
the next occurrence of *) into pure comment.

The compilers listed below have successfully compiled the program on
a 1 GHz Intel PIII running MS Windows 98 SE and partly on a 500 MHz AMD
K6-II running MS DOS 5.0; the respective execution times, roughly measured
with some overhead utilizing the PC clock (resolution: 0.055 s), are given.

Turbo Pascal 1.0 : 1.373 s / 2.527 s

Turbo Pascal 3.02 : 1.428 s / 2.636 s

Turbo Pascal 5.5  : 2.636 s / 2.856 s

Turbo Pascal 7.01 : 0.879 s / 0.989 s

Borland Delphi 4.0 build 5.37 : 0.824 s / -

Virtual Pascal 2.1 build 279 : 0.989 s / -

Free Pascal 2.4.4 i386-Win32 : 0.769 s / -

\vspace{0.1in}

To avoid a compilation result depending on the settings you could (not
in the case of Turbo Pascal versions 1.0 and 3.02) use the command line
versions of the compilers (TPC for Turbo Pascal, BPC for Borland Pascal 7,
DCC32 for Borland Delphi (32 bit versions; do not miss to use the -CC option
in order to generate a console executable), VPC for Virtual Pascal, FPC for
Free Pascal) instead of the compilers integrated in the IDEs.

\vspace{0.1in}

The program uses the heap and pointers in general just to store the vectors
constituting the adjacency matrix. It does it in a way reducing the readabilty
as little as possible while still allowing the use of compilers storing each
Boolean value in a separate byte and unable to handle memory block sizes
reaching or exceeding 64 KB.

If you don't change the respective compiler directives, range checks and (if
the compiler is compatible with Turbo Pascal 4.0) stack overflow checks are
generated. So the resulting executable will be very safe. It is also
quite small and needs less then 1 KB for the stack. Turbo Pascal versions
before 4.0 do not know conditional compiling and interpret some compiler
options differently. For enabling or disabling the generation of stack
overflow checks, the relevant letter is K (instead of S). So you could
insert the sequence \{\$K+\} (e.g. after \{\$R+\}) in order to generate those
checks.

\vspace{0.1in}

The program ignores any command line parameters or inputs other than
pressing Ctrl-C to cancel the execution.

It writes only to the standard output device. In the default case that
will be the monitor screen. But if the used compiler was compatible with
Turbo Pascal 4.0 you can redirect the output to a file.
That way the lines below enclosed in $<<<$ and $>>>$ were generated. You
could use them to compare your results with.

\vspace{0.1in}

$<<<$
\ttfamily
\small

=== Constructing and checking a G2(4) graph representation ===

===    Version 1   Copyright (c) 2013-07-31 Thomas Jenrich ===
\vspace{0.1in}

Search for orthogonal bases ... OK

Filling iso-sets ... OK

Allocating dynamic memory ... OK

Filling the adjacency matrix ... OK

Separating subsets B\_1, B\_2, B\_3, C ... OK

Checking the counts of adjacencies in B\_1, B\_2, B\_3 ... OK

OK

== Regular program stop ==

\rmfamily

$>>>$

\vspace{0.1in}

Author's eMail address: thomas.jenrich@gmx.de


\begin{thebibliography}{00}

\bibitem{Bor}
Karol Borsuk, \emph{Drei S\"{a}tze \"{u}ber die $n$-dimensionale euklidische
Sph\"{a}re}, Fund. Math., 20 (1933), 177-190.

\bibitem{Bon}

Andriy V. Bondarenko, \emph{On Borsuk's conjecture for two-distance sets},
arXiv:math.MG/1305.2584 (version 2),
http://front.math.ucdavis.edu/1305.2584v2

\bibitem{Bon2}

Andriy Bondarenko, \emph{On Borsuk's Conjecture for Two-Distance Sets},
Discrete Comput. Geom., Springer 2014, DOI 10.1007/s00454-014-9579-4

\bibitem{BroHae}

Andries E. Brouwer, Willem H. Haemers, \emph{Spectra of graphs},
 Monograph, 2012, Springer

\bibitem{Bro1}
Andries E. Brouwer, \emph{A construction of the Suzuki Graph}, 2008-07-08,

http://www.win.tue.nl/\verb+~+aeb/preprints/Suz.pdf

\bibitem{Bro2}

Andries E. Brouwer, \emph{G2(4) graph},

http://www.win.tue.nl/\verb+~+aeb/graphs/G24.html, retrieved 2013-05-22

\bibitem{HoKiNa}
N. Horiguchi, M. Kitazume, H. Nakasora,
\emph{A construction of the sporadic Suzuki graph from $U_3(4)$}, preprint,
2008

\bibitem{BrHoKiNa}
Andries E. Brouwer, Naoyuki Horiguchi, Masaaki Kitazume, Hiroyuki Nakasora,

\emph{A construction of the sporadic Suzuki graph from $U_3(4)$},

Journal of Combinatorial Theory, Series A 116 (2009), 1056-1062

\bibitem{Soi}
Leonard H. Soicher, \emph{ Three new distance-regular graphs},

Europ. J. Combin. 14 (1993), 501-505

\bibitem{CrnMik}
Dean Crnkovi\'{c}, Vedrana Mikuli\'{c},
\emph{BLOCK DESIGNS AND STRONGLY REGULAR GRAPHS CONSTRUCTED FROM THE
GROUP $U(3,4)$}, GLASNIK MATEMATI\v{C}KI, Vol. 41(61)(2006),

189-194

\bibitem{Che}
William Cherowitzo, \emph{A Short Introduction to Projective Geometry},

http://www-math.ucdenver.edu/\verb+~+wcherowi/courses/m5410/pgintro.pdf

\bibitem{Cam}
Robert Campbell, \emph{Lectures 12 \& 13: Finite Fields},
(Course) Number Theory, Math 413, 26 Jan 2003
http://userpages.umbc.edu/\verb+%7Ercampbel/Math413Spr05/Notes/12-13_Finite_Fields.html+

\bibitem{CibLip}
Namik Ciblak, Harvey Lipkin, \emph{ORTHONORMAL ISOTROPIC VECTOR BASES},
Proceedings of the 1998 ASME Design Engineering Technical Conference,
September, 13-16, 1998, Atlanta, Georgia

\bibitem{TP}\emph{Turbo Pascal} versions 1.0, 3.02, and 5.5 (binaries only)

http://edn.embarcadero.com/museum/antiquesoftware\#

For downloading one has to register or sign-in.

\bibitem{VP}\emph{Virtual Pascal} (Closed Source freeware)

One ZIP-file including binaries and documentation for Win32, OS/2, and Linux

Official forum:

http://vpascal.ning.com/

Forum entry \emph{Where can I download VP?} :

http://vpascal.ning.com/forum/topic/show?id=854411\verb+%3ATopic%3A9+

\bibitem{FP}\emph{Free Pascal} (Open Source freeware)

Sources, documentation, and binaries for several systems

http://www.freepascal.org

\end{thebibliography}
\end{document}